\documentclass{article}

\usepackage{amssymb}

\usepackage{graphicx}

\newtheorem{theorem}{Theorem}[section]

\newtheorem{lemma}[theorem]{Lemma}
\newtheorem{proposition}[theorem]{Proposition}

\begin{document}
	
\title{Solution manifolds of differential systems with discrete state-dependent delays are almost graphs} 
	
\author{Tibor Krisztin\thanks{This research was supported by the National Research, Development and Innovation Fund,
Hungary, grants NKFIH-K-129322, TKP2021-NVA-09.}
\\ELKH-SZTE Analysis and Applications Research Group\\ 
Bolyai Institute, University of Szeged \and Hans-Otto Walther\thanks{HOW is grateful for a stay at Bolyai Institute, University of Szeged, May 2022}\\Mathematisches Institut\\Universit\"at Gie{\ss}en}

\date{\small{\tt krisztin@math.u-szeged.hu}\\{\tt
		Hans-Otto.Walther@math.uni-giessen.de}}

\maketitle

\begin{abstract}
\noindent
We show that for a system 
$$
x'(t)=g(x(t-d_1(Lx_t)),\dots,x(t-d_k(Lx_t)))
$$
of $n$ differential equations with $k$ discrete state-dependent delays
the solution manifold, on which solution operators are differentiable, is nearly as simple as a graph over a closed subspace in $C^1([-r,0],\mathbb{R}^n)$. The map $L$ is continuous and linear from $C([-r,0],\mathbb{R}^n)$ onto a finite-dimensional vectorspace, and $g$ as well as the delay functions $d_{\kappa}$ are assumed to be continuously differentiable. 
\end{abstract}
	
\bigskip
	
	

\section{Introduction}

For an integer $n>0$ and a real number $r>0$ let $C_n=C([-r,0],\mathbb{R}^n)$ and $C^1_n=C^1([-r,0],\mathbb{R}^n)$
denote the Banach spaces of continuous and continuously differentiable maps $[-r,0]\to\mathbb{R}^n$, respectively, with the norms given by
$$
|\phi|_C=\max_{-r\le t\le0}|\phi(t)|\quad\mbox{and}\quad|\phi|=|\phi|_C+|\phi'|_C,
$$
for a chosen norm on $\mathbb{R}^n$. If a differential equation with state-dependent delay, like for example the equation
\begin{equation}
x'(t)=g(x(t-\Delta)),\quad\Delta=d(x(t))
\end{equation}
 with real functions $g:\mathbb{R}\to\mathbb{R}$ and $d:\mathbb{R}\to[0,r]$,
is written in the general form
\begin{equation}\label{gen-equ}
x'(t)=f(x_t)
\end{equation}
of an autonomous delay differential equation, with a map $f: U\to\mathbb{R}^n$ defined on an open subset of $C^1_n$, and  $x_t(s)=x(t+s)$, $-r\le s\le 0$, then the associated {\it solution manifold} is the set 
$$
X_f=\{\phi\in U:\phi'(0)=f(\phi)\}.
$$ 
If $f$ is continuously differentiable and satisfies the extension property that

\medskip

(e) {\it each derivative $Df(\phi):C^1_n\to\mathbb{R}^n$, 
$\phi\in U$, continues to a linear map 
$D_ef(\phi):C_n\to\mathbb{R}^n$ and the map 
$$
U\times C_n\ni(\phi,\chi)\mapsto D_ef(\phi)\chi)\in\mathbb{R}^n
$$
is continuous}

\medskip

and if $X_f$ is non-empty then $X_f$ is a continuously differentiable submanifold of codimension $n$ in $C^1_n$, and the initial value problem for initial data in $X_f$ is well-posed, with each solution operator continuously differentiable \cite{W1,HKWW}. 

\medskip

The extension property (e) is a version of the notion of being {\it almost Fr\'echet differentiable} which was introduced for maps on domains $U\subset C_n$ by Mallet-Paret, Nussbaum, and Paraskevopoulos \cite{M-PNP}.

\medskip

Let us recall that in contrast to theory 
for differential equations with constant time lags \cite{HVL,DvGVLW} the initial value problem for equations with state-dependent delays is in general not well-posed for initial data in open subsets of the space $C_n$.

\medskip

The present paper about solution manifolds continues work which started with \cite{KR,W6}. Under a boundedness condition on the extended derivatives $D_ef$, or under a condition on $f$ which generalizes delays being bounded away from zero, solution manifolds are graphs which can be written as
$$
X_f=\{\chi+\alpha(\chi):\chi\in dom\}
$$
with a continuously differentiable map $\alpha$ from an open subset $dom$ of the closed subspace
$$
X_0=\{\phi\in C^1_n:\phi'(0)=0\}
$$
into a complementary space $Q\subset C^1_n$, see the proof of
\cite[Lemma 1]{KR} and \cite[Theorem 2.4]{W6}, respectively. 

\medskip 

An example of the form (1) with $d:\mathbb{R}\to(0,r]$ in \cite[Section 3]{W6} shows that in general solution manifolds do not admit any graph representation. However, \cite[Theorem 5.1]{W6} says that for a class of systems with discrete state-dependent delays all of which are strictly positive the solution manifolds are nearly as simple as a graph, namely, they are {\it almost graphs over} $X_0$ in the sense of the following definition from \cite[Section 1]{W7}: 

\medskip

{\it A continuously differentiable submanifold $X$ of a Banach space $E$ is called an almost graph over a closed subspace $H\subset E$ if $H$ has a closed complementary subspace in $E$ and if there is a continuously differentiable map $\alpha:H\supset dom\to E$, $dom$ an open subset of $H$, such that
\begin{eqnarray}
X & = & \{\zeta+\alpha(\zeta)\in E:\zeta\in dom\},\nonumber\\
\alpha(\zeta) & = & 0\quad\mbox{on}\quad dom\cap X,\nonumber\\ \alpha(\zeta) & \in & E\setminus H\quad\mbox{on}\quad dom\setminus X,\nonumber
\end{eqnarray}
and the map $H\supset dom\ni\zeta\mapsto\zeta+\alpha(\zeta)\in X$ is  a diffeomorphism onto $X$.}

\medskip

An example of an almost graph in the plane is the unit circle without $(0,1)^{tr}$, for $H=\mathbb{R}\times\{0\}$ and $\alpha$ the inverse of the stereographic projection \cite[Section 1]{W7}. 

\medskip

A property stronger than being an almost graph is existence of an {\it almost graph diffeomorphism} as introduced in \cite[Section 1]{W7}: For $E,X,H$ as above, and for an
open set ${\mathcal O}\subset E$,

\medskip

{\it  a diffeomorphism $A:{\mathcal O}\to E$ onto an open subset of $E$ is called an almost graph diffeomorphism with respect to $X$, $H$, and ${\mathcal O}$, if}\footnote{as in the definition of a continuously differentiable submanifold}
$$
A(X\cap {\mathcal O})=H\cap A({\mathcal O})
$$
{\it and if  $A$ leaves the points of $(X\cap {\mathcal O})\cap H$ fixed.}

\medskip 

In \cite[Section 1]{W7} it is shown that the existence of an almost graph diffeomorphism with respect to $X,H,{\mathcal O}$ implies
that $X\cap {\mathcal O}$ is an almost graph over $H$.

\medskip

The main results of \cite{W7}, Theorems 3.5 and 4.8, imply that for a
class of systems more general than those studied in \cite{W6} the associated solution manifolds carry a finite atlas of manifold charts
whose domains are almost graphs over $X_0$. The size of the atlas is determined precisely by the zerosets of the delays. If for example the delay function $d$ in Eq. (1) is non-constant and has zeros then the atlas found in \cite{W7} consists of exactly 2 manifold charts.

\medskip

Let us now recall the systems
\begin{equation}\label{sdde}
x'(t)=g(x(t-d_1(Lx_t)),\ldots,x(t-d_{k}(Lx_t)))
\end{equation}
introduced in \cite{W7}. In Eq. (3) the delays are given by compositions of a continuous linear map $L:C_n\to F$ onto a 
finite-dimensional normed vectorspace over the field $\mathbb{R}$,
with continuously differentiable delay functions $d_{\kappa}:W\to[0,r]$, $W\subset F$ open and $\kappa\in\{1,\ldots,k\}$. The nonlinearity $g$ is a continuously differentiable map from an open subset $V\subset\mathbb{R}^{nk}$ into $\mathbb{R}^n$.
 
\medskip

The notation for the argument of $g$ in Eq. (3) is an abbreviation for the column vector $y\in\mathbb{R}^{nk}$ with components
$$
y_{\iota}=x_{\nu}(t-d_{\kappa}(Lx_t))\quad(\quad=(x_t)_{\nu}(-d_{\kappa}(Lx_t))\quad)
$$
for $\iota=(\kappa-1)n+\nu$ with $\kappa\in\{1,\ldots,k\}$ and $\nu\in\{1,\ldots,n\}$. 

\medskip

With regard to the form of the delays in Eq. (3) one may think of $L\phi$ as an approximation of $\phi\in C_n$ in the subspace $F\subset C_n$, and view $d_{\kappa}(L\cdot)$ as a substitute for a more general delay functional defined on an open subset of $C_n$.

\medskip

In order that Eq. (3) makes sense it is assumed that 

\medskip

(V) {\it there exist $\phi\in C^1_n$ with $L\phi\in W$ and
	$(\phi(-d_1(L\phi)),\ldots,\phi(-d_k(L\phi)))\in V$.}

\medskip

(in notation as described above for $\phi=x_t$).

With 
$$
\widehat{\phi}=(\phi(-d_1(L\phi)),\ldots,\phi(-d_k(L\phi)))
$$
we get that
$$
U=\{\phi\in C^1_n:L\phi\in W\quad\mbox{and}\quad\widehat{\phi}\in V\}
$$
is non-empty, and Eq. (3) takes the form of Eq. (2)
with $f:U\to\mathbb{R}^n$ given by
$$
f(\phi)=g(\widehat{\phi}).
$$
According to \cite[Propositions 2.1 and 2.3]{W7} the set $U$ is open, $f$ is continuously differentiable with property (e), and $X_f\neq\emptyset$, so that $X_f$ is a continuously differentiable submanifold of codimension $n$ in $C^1_n$.

\medskip

In the present paper we consider system (\ref{sdde}) under the above conditions and construct an almost graph diffeomorphism with respect to the whole solution manifold $X_f$, to the subspace $X_0$, and to the open neighbourhood ${\mathcal O}=U$ of $X_f$ in $C^1_n$. The result is stated as Theorem 3.5 below. 

\medskip

Let us emphasize that the finite atlas result from \cite{W7} holds true under further hypotheses on $g$ or on the delay functions $d_{\kappa}$ whereas the construction in Sections 2-3 below
requires nothing beyond smoothness as stated above. This discrepancy reflects the fact that in the proof
of \cite[Theorem 5.1]{W7}, as well as in the proofs of \cite[Theorems 2.4 and 5.1]{W6}, the invariance property 
$$
\widehat{A(\phi)}=\widehat{\phi}
$$
of a diffeomorphism $A$ from $X_f$ onto an open subset of $X_0$ is established and exploited. The approach in the present paper, which originated in the case study  \cite{W9}, proceeds without recourse to the said invariance property.

\medskip 

The Introduction of \cite{W7} lists a wide variety of special cases of system (\ref{sdde}). Therefore, by Theorem 3.5, the solution manifolds associated with many familiar delay differential systems with state-dependent delays are almost graphs over a closed subspace of  $C^1_n$. 
It is an open problem whether this is true for the solution manifold of the general equation (\ref{gen-equ}) with a continuously differentiable $f:U\to \mathbb{R}^n$, defined on an open $U\subset C^1_n$, satisfying the extension property (e). 

\medskip 

There are differential equations with discrete delays so that the delay functions are not of the form $\delta_\kappa (L\phi)$ with a continuous linear map 
$L:C_n\to F$ into a finite dimensional vectorspace $F$. 
For example, threshold delays and transmission delays are implicitly defined (see e.g. \cite{A,KS1,KS2,W2,W3,HKWW,KW}), and the corresponding delay functions $\sigma(\phi)$ are defined on an open subset of $C^1_n$ in order to write the system in the form (\ref{gen-equ}) with the required properties. 
It would be interesting to describe the solution manifolds for systems with threshold and transmission delays as well.

\medskip

{\bf Notation, conventions, preliminaries.} For subsets $A\subset B$ of a topological space $T$ we say $A$ is open in $B$ if $A$ is open with respect to the relative topology on $B$. Analogously for $A$ closed in $B$. The relation $A\subset\subset B$ for open subsets of $T$  means that the closure $\overline{A}$ of $A$ is compact and contained in $B$.

\medskip

Finite-dimensional vectorspaces are always equipped with the canonical topology which makes them topological vectorspaces. 

\medskip

On $\mathbb{R}^{nk}$ we use a norm which satisfies
$$
|y_j|\le|y|\le\sum_{\iota=1}^{nk}|y_{\iota}|\quad\mbox{for all }\quad y\in\mathbb{R}^{nk},\,\,j\in\{1,\ldots,nk\}.
$$

\medskip

In case $V\neq\mathbb{R}^{nk}$ the expression 
$dist(v,\mathbb{R}^{nk}\setminus V)=\min_{y\in\mathbb{R}^{nk}\setminus V}|y-v|$ defines a continuous function $V\to(0,\infty)$.

\medskip

An upper index as in $(x_1,\ldots,x_N)^{tr}\in\mathbb{R}^N$ denotes the transpose of the row vector $(x_1,\ldots,x_N)$. Vectors in $\mathbb{R}^N$ which occur as argument of a map are always written as row vectors. The vectors of the canonical basis of $\mathbb{R}^N$ are denoted by $e_{\nu}$, $\nu\in\{1,\ldots,N\}$; $e_{\nu\mu}=1$ for $\nu=\mu$ and $e_{\nu\mu}=0$ for $\nu\neq\mu$, $\nu$ and $\mu$ in $\{1,\ldots,N\}$.   

\medskip

Derivatives  and partial derivatives of a map at a given argument are continuous linear maps, indicated by a capital $D$. In case of real functions on domains in $\mathbb{R}$ and in $\mathbb{R}^N$, $\phi'(t)=D\phi(t)1$ and $\partial_{\nu}g(x)=D_{\nu}g(x)1$, respectively.

\medskip

For $n=1$ we abbreviate $C=C_1$ and $C^1=C^1_1$.

\medskip

We define continuous bilinear products $C\times\mathbb{R}^n\to C_n$ 
and $\mathbb{R}^n\times C_n\to C_n$ by
$$
(\phi\cdot q)_{\nu}=q_{\nu}\phi\in C
$$ 
and
$$
(q\cdot\phi)_{\nu}=q_{\nu}\phi_{\nu}\in C
$$ 
for $\nu=1,\ldots,n$. Obviously
$$
q\cdot\phi=\sum_{\nu=1}^nq_{\nu}(\phi_{\nu}\cdot e_{\nu})
$$
for all $q\in\mathbb{R}^n,\phi\in C_n$.

\medskip



The maps
$$
C^1_n\times W\ni(\psi,\eta)\mapsto\psi_{\nu}(-d_{\kappa}(\eta))\in\mathbb{R},\quad\nu\in\{1,\ldots,n\},\quad\kappa\in\{1,\ldots,k\},
$$
are continuously differentiable, compare Part 2.1 of the proof of \cite[Proposition 2.1]{W7}. It follows that also the map
$$
U\ni\phi\mapsto\widehat{\phi}\in V\subset\mathbb{R}^{nk}
$$
is continuously differentiable.

\medskip

The inclusion $C^1_n\ni\phi\mapsto\phi\in C_n$,
differentiation $\partial:C^1_n\ni\phi\mapsto \phi'\in C_n$, and  
evaluation  
$ev_t:C\ni\phi\mapsto\phi(t)\in\mathbb{R}$ at $t\in[-r,0]$ are  continuous linear maps.

\medskip


\medskip


In the sequel a diffeomorphism is a continuously differentiable injective map with open image whose inverse is continuously differentiable.

\medskip

\section{Preparations}

The following lemma is a version of \cite[Proposition 4.1]{W6} which includes smallness in $C$ of the function found.

\begin{lemma} 
Let $\lambda:C\to{\mathcal F}$ be a continuous linear map into a finite-dimensional real vectorspace ${\mathcal F}$, and let $\epsilon>0$ be given. There exists $\psi\in C^1$ with $\lambda\psi=0$, $\psi'(0)=1$, and $|\psi|_C<\epsilon$.
\end{lemma}

{\bf Proof.} 1. Proof that there exists a complementary space $K\subset C^1$ for  $\lambda^{-1}(0)$ in $C$. 
We have
\begin{eqnarray*}
	\lambda C^1 & = & \overline{\lambda C^1}\quad\mbox{(with}\quad\,\dim\,\lambda C^1\le\dim\,\lambda C<\infty)\\
	& \supset & \lambda C\quad\mbox{(with}\quad  C^1\quad\mbox{dense in}\quad C)\\
	& \supset  & \lambda C^1,
\end{eqnarray*} 
hence $\lambda C^1=\lambda C$. Set $K=\sum_{j=1}^k\mathbb{R}\psi_j$ with preimages $\psi_j\in C^1$ of a basis of $\lambda C$ and verify $C=\lambda^{-1}(0)\oplus K$.

\medskip

2.  The projection $P:C\to C$ along $K\subset C^1$ onto $\lambda^{-1}(0)$ maps
$C^1$ into $\lambda^{-1}(0)\cap C^1$. Choose a sequence $(\phi_m)_1^{\infty}$ in $C^1$ with $\phi_m'(0)=1$ for all $m\in\mathbb{N}$ and $|\phi_m|_C\to0$ as $m\to\infty$. Then $|(id-P)\phi_m|_C\to0$ as $m\to\infty$. As $K$ is finite-dimensional we also get $|(id-P)\phi_m|\to0$ as $m\to\infty$. It follows that $ev_0\partial(id-P)\phi_m\to0$ as $m\to\infty$.
Using this and $P\phi_m\in C^1$ we get
$$
1=\phi_m'(0)=ev_0\partial\phi_m=ev_0\partial(P\phi_m+(id-P)\phi_m)=ev_0\partial P\phi_m+ev_0\partial(id-P)\phi_m
$$
for all $m\in\mathbb{N}$, which yields $ev_0\partial P\phi_m\to1$ as $m\to\infty$. For $m$ so large that $ev_0\partial P\phi_m\neq0$
the functions
$$
\psi_m=\frac{1}{ev_0\partial P\phi_m} P\phi_m
$$ 
belong to $C^1$ and satisfy $\lambda\psi_m=0$ and $\psi_m'(0)=ev_0\partial\psi_m=1$. For $m$ sufficiently large we also obtain $|\psi_m|_C<\epsilon$. $\Box$

\medskip

The next proposition provides functions in $C^1$ which after multiplication with the unit vectors $e_{\nu}\in\mathbb{R}^n$ yield bases of subspaces which are complementary for $X_0$ and depend on $\phi\in U$ via $v=\widehat{\phi}\in V\subset\mathbb{R}^{nk}$. For $\phi\in X_f$ the subspace is complementary also for the tangent space $T_{\phi}X_f$. We omit proofs of these facts as they will not be used in the sequel. The almost graph diffeomorphism associated with $X_f$, $X_0$, and $U$ which will be constructed in the next section is composed of translations by vectors in the complementary spaces just mentioned. 

\medskip

Recall the map $L:C_n\to F$ from Eq. (3).

\begin{proposition}
Let a continuous function $h:V\to(0,\infty)$ and $\nu\in\{1,\ldots,n\}$ be given. Then there exists a continuously differentiable map
$$
H_{\nu}:V\to C^1
$$
so that for all $v\in V$,
$$
L(H_{\nu}(v)\cdot e_{\nu})=0,\,\,(H_{\nu}(v))'(0)=1,\,\,|H_{\nu}(v)|_C\le h(v),
$$ 
and for each $\mu\in\{1,\ldots,nk\}$,
$$
|D_{\mu}H_{\nu}(v)1|_C\le h(v).
$$
\end{proposition}

{\bf Proof.} 1. There is a sequence of non-empty open subsets $V_{j1},V_{j2},V_j$ of $V$, with $j\in\mathbb{N}$, such that
$$
\bigcup_{j=1}^{\infty}V_j=V,
$$
and for every $j\in\mathbb{N}$,
$$
V_{j1}\subset\subset V_{j2}\subset\subset V_j\quad\mbox{and}\quad V_j\subset\subset V_{j+1,1}.
$$
With $\overline{V_{02}}=\emptyset$ we have that for each integer $j\ge 1$,
$$
(V_{j+1}\setminus\overline{V_{j2}})\cap(V_j\setminus\overline{V_{j-1,2}})=V_j\setminus\overline{V_{j2}}
$$
while for integers $j\ge1$ and $k\ge j+2$, 
$$
(V_k\setminus\overline{V_{k-1,2}})\cap(V_j\setminus\overline{V_{j-1,2}})=V_j\setminus\overline{V_{k-1,2}}\subset V_j\setminus \overline{V_{j+1,2}}=\emptyset. 
$$
 
2. For every $j\in\mathbb{N}$ choose a continuously differentiable function
$$
a_j:\mathbb{R}^{nk}\to[0,1]
$$
with
$$
a_j(v)=1\,\,\mbox{on}\,\,\overline{V_{j1}},\quad a_j(v)=0\,\,\mbox{on}\,\,\mathbb{R}^{nk}\setminus V_{j2}.
$$
For  every $j\in\mathbb{N}$ choose an upper bound
$$
A_j>1+\sum_{\mu=1}^{nk}\max_{v\in\mathbb{R}^{nk}}|D_{\mu}a_j(v)1|=\max_{v\in\overline{V_j}}|a_j(v)|+\sum_{\mu=1}^{nk}\max_{v\in\overline{V_j}}|D_{\mu}a_j(v)1|
$$
so that the sequence $(A_j)_1^{\infty}$ in $[1,\infty)$ is increasing.

\medskip

The sequence $(h_j)_1^{\infty}$ given by $h_j=\min_{v\in\overline{V_j}}h(v)>0$ is nonincreasing. We have 
$$
\frac{h_j}{2A_j}\le h_j\,\,\mbox{for all}\,\,j\in\mathbb{N},
$$
and the sequence $(h_j/2A_j)_{j=1}^{\infty}$ is decreasing.

\medskip

3. For each $j\in\mathbb{N}$ apply  Lemma 2.1 to $\lambda:C\to F$ given by $\lambda\phi=L(\phi\cdot e_{\nu})$, and to $\epsilon=h_j/2A_j$. This yields a sequence of functions $\psi_{\nu,j}\in C^1$, $j\in\mathbb{N}$,  which satisfy
$$
L(\psi_{\nu,j}\cdot e_{\nu})=0,\,\,(\psi_{\nu,j})'(0)=1,\,\,|\psi_{\nu,j}|_C<\frac{h_j}{2A_j}.
$$
The maps
$$
H_{\nu,j}:V_j\setminus\overline{V_{j-1,2}}\to C^1,\quad j\in\mathbb{N},
$$
given by
$$
H_{\nu,j}(v)=a_j(v)\psi_{\nu,j}+(1-a_j(v))\psi_{\nu,j+1}
$$
are continuously differentiable and satisfy
\begin{eqnarray}
L(H_{\nu,j}(v)\cdot e_{\nu}) & = & a_j(v)L(\psi_{\nu,j}\cdot e_{\nu})+(1-a_j(v))L(\psi_{\nu,j+1}\cdot e_{\nu})=0,\nonumber\\
(H_{\nu,j}(v))'(0) & = & a_j(v)\psi_{\nu,j}'(0)+(1-a_j(v))\psi_{\nu,j+1}'(0)=1,\nonumber\\
|H_{\nu,j}(v)|_C & \le & a_j(v)|\psi_{\nu,j}|_C+(1-a_j(v))|\psi_{\nu,j+1}|_C\nonumber\\
& \le & a_j(v)h_j+(1-a_j(v))h_{j+1}\nonumber\\
& \le & a_j(v)h_j+(1-a_j(v))h_j=h_j\le h(v)\nonumber
\end{eqnarray} 
for every $j\in\mathbb{N}$ and all $v\in V_j\setminus\overline{V_{j-1,2}}$. Moreover, for such $j$ and $v$, and for every $\mu\in\{1,\ldots,kn\}$,
$$
D_{\mu}H_{\nu,j}(v)1=(D_{\mu}a_j(v)1)\,\psi_{\nu,j}- (D_{\mu}a_j(v)1)\,\psi_{\nu,j+1},
$$
hence
\begin{eqnarray}
|D_{\mu}H_{\nu,j}(v)1|_C & \le & |D_{\mu}a_j(v)1|(|\psi_{\nu,j}|_C+|\psi_{\nu,j+1}|_C)\nonumber\\
& \le & A_j\left(\frac{h_j}{2A_j}+\frac{h_{j+1}}{2A_{j+1}}\right)\le A_j\left(2\frac{h_j}{2A_j}\right)=h_j\le h(v).\nonumber
\end{eqnarray}

\medskip

4. It remains to show that the maps $H_{\nu,j}$, $j\in\mathbb{N}$, define a map $H_{\nu}:V\to C^1$. This follows from Part 1 of the proof provided $H_{\nu,j+1}$ and $H_{\nu,j}$ coincide on the intersection $V_j\setminus\overline{V_{j2}}$ of their domains, for
every $j\in\mathbb{N}$. For $j\in\mathbb{N}$ and $v\in V_j\setminus\overline{V_{j2}}$ we have
$$
H_{\nu,j+1}(v)=a_{j+1}(v)\psi_{\nu,j+1}+(1-a_{j+1}(v))\psi_{\nu,j+2}=\psi_{\nu,j+1}
$$
due to $a_{j+1}(v)=1$ on $\overline{V_{j+1,1}}\supset V_j\supset V_j\setminus\overline{V_{j2}}$,
and
$$
H_{\nu,j}(v)=a_j(v)\psi_{\nu,j}+(1-a_j(v))\psi_{\nu,j+1}=\psi_{\nu,j+1}
$$
due to $a_j(v)=0$ on $\mathbb{R}^{nk}\setminus\overline{V_{j2}}\supset V_j\setminus\overline{V_{j2}}$. $\Box$

\medskip

\section{The almost graph diffeomorphism}

\medskip

The function $h_V:V\to(0,\infty)$ given by
$$
h(v)=\frac{\min\{1,dist(v,\mathbb{R}^{nk}\setminus V)\}}{2(nk)^2(1+\max_{\iota=1,\ldots,nk;\nu=1,\ldots,n}|\partial_{\iota}g_{\nu}(v)|+\max_{\nu=1,\ldots,n}|g_{\nu}(v)|)}.
$$
in case $V\neq\mathbb{R}^{nk}$ and 
$$
h(v)=\frac{1}{2(nk)^2(1+\max_{\iota=1,\ldots,nk;\nu=1,\ldots,n}|\partial_{\iota}g_{\nu}(v)|+\max_{\nu=1,\ldots,n}|g_{\nu}(v)|)}
$$
for $V=\mathbb{R}^{nk}$ is continuous.

\medskip

For $h=h_V$ choose functions $H_{\nu}:V\to C^1$, $\nu\in\{1,\ldots,n\}$,  according to Proposition 2.2 and define $H:V\to C^1_n$ by 
$H(v)=\sum_{\nu=1}^nH_{\nu}(v)\cdot\,e_{\nu}$, or equivalently,
$(H(v))_{\nu}=H_{\nu}(v)$ for $\nu=1,\ldots,n$. The map $A:U\to C^1_n$
given by
$$
A(\phi)=\phi-g(v)\cdot H(v)=\phi-\sum_{\nu=1}^ng_{\nu}(v)(H_{\nu}(v)\cdot\,e_{\nu})\quad\mbox{with}\quad v=\widehat{\phi}
$$
is continuously differentiable and satisfies $A(X_f)\subset X_0$ as for every $\phi\in X_f$ we have 
\begin{eqnarray}
(A(\phi))'(0) & = & \phi'(0)-(g(v)\cdot H(v))'(0)\quad\mbox{(with}\quad v=\widehat{\phi})\nonumber\\
& =  & \phi'(0)-\left(\sum_{\nu=1}^ng_{\nu}(v)(H_{\nu}(v)\cdot\,e_{\nu})\right)'(0)\nonumber\\
& =  & \phi'(0)-\sum_{\nu=1}^ng_{\nu}(v)\,(H_{\nu}(v))'(0)\,e_{\nu}\nonumber\\
& = & \phi'(0)-g(v)\nonumber\\
& = & \phi'(0)-g(\widehat{\phi})=0\quad\mbox{(with}\quad\phi\in X_f),\nonumber
\end{eqnarray}
which means $A(\phi)\in X_0$. 
From the above lines it is also obtained that if $\phi\in U$ and $A(\phi)\in X_0$, then   $0=(A(\phi))'(0)=\phi'(0)-g(\widehat{\phi})$, that is $\phi\in X_f$. 

\medskip

We also have $A(\phi)=\phi$ on $X_0\cap X_f$ since $\phi\in X_0\cap X_f$ yields $0=\phi'(0)=g(\widehat{\phi})$, hence $A(\phi)=\phi-g(\widehat{\phi})\cdot H(\widehat{\phi})=\phi$.

\medskip

For $A$ to be an almost graph diffeomorphism associated with the submanifold $X_f\subset C^1_n$, with the open set $U\subset C^1_n$, and with the closed subspace $X_0\subset C^1_n$, it remains to prove that $A$ is a diffeomorphism onto an open subset of $C^1_n$.

\medskip

Observe that due to $L(H_{\nu}(v)\cdot\,e_{\nu})=0$ we have
\begin{equation}
LA(\phi)=L\phi-L\left(\sum_{\nu=1}^ng_{\nu}(v)(H_{\nu}(v)\cdot\,e_{\nu})\right)=L\phi
\end{equation}
for every $\phi\in U$.

\medskip

Next we examine the relation between $v=\widehat{\phi}$, $\phi\in U$, and $y=\widehat{\chi}$ for $\chi=A(\phi)$. Let
$\eta=L\chi=L\phi\in W$. For $\iota=(\kappa-1)n+\nu  $ with $\kappa\in\{1,\ldots,k\}$ and $\nu\in\{1,\ldots,n\}$,
\begin{eqnarray}
y_{\iota} & = & \widehat{\chi}_{\iota}=\widehat{A(\phi)}_{\iota}= [\phi_{\nu}-(g(v)\cdot H(v))_{\nu}](-d_{\kappa}(\eta))\nonumber\\
& = & \widehat{\phi}_{\iota}-(g_{\nu}(v)(H_{\nu}(v))(-d_{\kappa}(\eta))\nonumber\\
& = & v_{\iota}-g_{\nu}(v)(H_{\nu}(v))(-d_{\kappa}(\eta))\nonumber\\
& = & S_{\iota}(\eta,v)\nonumber
\end{eqnarray}
with the continuously differentiable map $S:W\times V\to\mathbb{R}^{nk}$ given by
$$
S(\eta,v)=v-R(\eta,v)
$$
and
\begin{eqnarray}
R_{\iota}(\eta,v) & = & g_{\nu}(v)(H_{\nu}(v))(-d_{\kappa}(\eta))\nonumber\\
& = & ev_{-d_{\kappa}(\eta)}(g_{\nu}(v)H_{\nu}(v))\nonumber
\end{eqnarray}
for $\iota=(\kappa-1)n+\nu$ with $\kappa\in\{1,\ldots,k\}$ and $\nu\in\{1,\ldots,n\}$.

\medskip

\begin{proposition}
(i) For all $(\eta,v)\in W\times V$, 
$$
|D_2R(\eta,v)|_{L_c(\mathbb{R}^{nk},\mathbb{R}^{nk})}\le\frac{1}{2}.
$$
(ii) In case $V\neq\mathbb{R}^{nk}$,
$$
|R(\eta,v)<\frac{dist(v,\mathbb{R}^{nk}\setminus V)}{2}
$$
for all $(\eta,v)\in W\times V$.
\end{proposition}

{\bf Proof.} 1. On assertion (i). Let $\eta\in W$ be given and define $R_{\eta}:V\to\mathbb{R}^{nk}$ by $R_{\eta}(v)=R(\eta,v)$. 
For every $v\in V$ and for all $y\in\mathbb{R}^{nk}$ with $|y|\le1$ we get
$$
|D_2R(\eta,v)y|=|DR_{\eta}(v)y|\le\sum_{\iota=1}^{nk}\sum_{j=1}^{nk}|\partial_jR_{\eta,\iota}(v)\cdot y_j|\le\sum_{\iota=1}^{nk}\sum_{j=1}^{nk}|\partial_jR_{\eta,\iota}(v)|,
$$
and for $j\in\{1,\ldots,nk\}$ and $\iota=(\kappa-1)n+\nu$ with $\kappa\in\{1,\ldots,k\}$ and $\nu\in\{1,\ldots,n\}$, by the chain rule,
\begin{eqnarray}
|\partial_jR_{\eta,\iota}(v)| & = & |D_jR_{\eta,\iota}(v)1|\nonumber\\
& = & |ev_{-d_{\kappa}(\eta)}\left(D_jg_{\nu}(v)(1)\cdot H_{\nu}(v)+g_{\nu}(v)\cdot D_jH_{\nu}(v)1\right)|\nonumber\\
& = & |D_jg_{\nu}(v)(1)\cdot(H_{\nu}(v))(-d_{\kappa}(\eta))+g_{\nu}(v)\cdot(D_jH_{\nu}(v)1)(-d_{\kappa}(\eta))|\nonumber\\
& \le & (|\partial_jg_{\nu}(v)|+|g_{\nu}(v)|)h(v)\quad\mbox{(with Proposition 2.2)}.\nonumber
\end{eqnarray}
It follows that
\begin{eqnarray}
|D_2R(\eta,v)|_{L_c(\mathbb{R}^{nk},\mathbb{R}^{nk})} & = & \sup_{|y|\le1}|D_2R(\eta,v)y|\nonumber\\
& \le & (nk)^2[\max_{j=1,\ldots,nk;\nu=1,\ldots,n}|\partial_jg_{\nu}(v)|+\max_{\nu=1,\ldots,n}|g_{\nu}(v)|]h(v)\nonumber\\
& \le & \frac{1}{2}\quad\mbox{(see the choice of}\,\,h).\nonumber
\end{eqnarray}

2. On assertion (ii). Assume $V\neq\mathbb{R}^{nk}$. For all $(\eta,v)\in W\times V$ and for each $\iota=(\kappa-1)n+\nu$ with
$\kappa\in\{1,\ldots,k\}$ and $\nu\in\{1,\ldots,n\}$ we have
\begin{eqnarray}
|R_{\iota}(\eta,v)| & = & |g_{\nu}(v)(H_{\nu}(v))(-d_{\kappa}(\eta))|\nonumber\\
& \le & |g_{\nu}(v)||H_{\nu}(v)|_C\le  |g_{\nu}(v)|h(v).\nonumber
\end{eqnarray}
Hence
\begin{eqnarray}
|R(\eta,v)| & \le & \sum_{\iota=1}^{nk}|R_{\iota}(\eta,v)|\le h(v)\sum_{\iota=1}^{nk}\max_{\nu=1,\ldots,n}|g_{\nu}(v)|\nonumber\\
& \le & h(v)\,nk\,\max_{\nu=1,\ldots,n}|g_{\nu}(v)|\nonumber\\
& < & \frac{dist(v,\mathbb{R}^{nk}\setminus V)}{2}.\quad\Box\nonumber
\end{eqnarray}

\medskip

It is convenient to introduce the continuously differentiable maps $S_{\eta}:V\ni v\mapsto S(\eta,v)\in\mathbb{R}^{nk}$, $\eta\in W$.

\medskip

\begin{proposition}
(i) The set $\cup_{\eta\in W}\{\eta\}\times S_{\eta}(V)\subset F\times\mathbb{R}^{nk}$ is open.

\medskip

(ii) Each map $S_{\eta}:V\to\mathbb{R}^{nk}$, $\eta\in W$, is a diffeomorphism onto the open set  $S_{\eta}(V)=S(\{\eta\}\times V)\subset\mathbb{R}^{nk}$.

\medskip

(iii) The map $\cup_{\eta\in W}\{\eta\}\times S_{\eta}(V)\ni(\eta,y)\mapsto S_{\eta}^{-1}(y)\in V$ is continuously differentiable.
\end{proposition}

\medskip

{\bf Proof.} 1. On assertion (i). Let $\eta_0\in W$ and $y_0=S_{\eta_0}(v_0)=S(\eta_0,v_0)$ with $v_0\in V$ be given. Choose a closed ball $V_0\subset V$ with center $v_0$ and radius $\epsilon>0$, and an open neighbourhood $W_0$ of $\eta_0$ in $W$ such that for all $\eta\in W_0$,
$$
|R(\eta,v_0)-R(\eta_0,v_0)|<\frac{\epsilon}{8}.
$$
For $\eta\in W_0$ and $v,v_1$ in $V_0$ the Mean Value Theorem in combination with Proposition 3.1 yield
$$
|R(\eta,v)-R(\eta,v_1)|\le\frac{1}{2}|v-v_1|.
$$
Let $Y\subset\mathbb{R}^{nk}$ denote the open ball with center $y_0$ and radius $\epsilon/8$. Let $\eta\in W_0$ and $y\in Y$ be given and consider the map
$$
V_0\ni v\mapsto y+R(\eta,v)\in\mathbb{R}^{nk},
$$
which is a contraction. For each $v\in V_0$,
\begin{eqnarray}
|y+R(\eta,v)-v_0| & = & |y+R(\eta,v)-(y_0+R(\eta_0,v_0)|\nonumber\\
 & \le & |y-y_0|+|R(\eta,v)-R(\eta,v_0)|+|R(\eta,v_0)-R(\eta_0,v_0)|\nonumber\\
& \le & \frac{\epsilon}{8}+\frac{1}{2}|v-v_0|+\frac{\epsilon}{8}\le\frac{3\epsilon}{4}\le\epsilon,\nonumber
\end{eqnarray}
and we see that the previous contraction has range in $V_0$. Consequently there is a fixed point
$$
v=y+R(\eta,v)\in V_0\subset V.
$$
Hence $y=v-R(\eta,v)=S_{\eta}(v)$. It follows that
$$
W_0\times Y\subset\cup_{\eta\in W}\{\eta\}\times S_{\eta}(V),
$$
which yields the assertion.

\medskip

2. On assertion (ii). Let $\eta\in W$ be given. It follows from assertion (i) that the set $S_{\eta}(V)\subset\mathbb{R}^{nk}$ is open. 

\medskip

2.1. Proof that $S_{\eta}$ is injective in case $V\neq\mathbb{R}^{nk}$. Let $v,\tilde{v}$ in $V$ be given with $S_{\eta}(v)=S_{\eta}(\tilde{v})$. Then
$$
v-\tilde{v}=R(\eta,v)-R(\eta,\tilde{v}).
$$
Without loss of generality, $dist(\tilde{v},\mathbb{R}^{nk}\setminus V)\le dist(v,\mathbb{R}^{nk}\setminus V)$. Using Proposition 3.1 (ii) we infer
$$
|\tilde{v}-v|=|R(\eta,\tilde{v})-R(\eta,v)|<\frac{1}{2}(dist(\tilde{v},\mathbb{R}^{nk}\setminus V)+dist(v,\mathbb{R}^{nk}\setminus V))
\le dist(v,\mathbb{R}^{nk}\setminus V).
$$
It follows that the line segment $v+[0,1](\tilde{v}-v)$ belongs to $V$, and the Mean Value Theorem in combination with Proposition 3.1 (i) yields
$$
|v-\tilde{v}|=|R(\eta,v)-R(\eta,\tilde{v})|\le\frac{1}{2}|v-\tilde{v}|,
$$
which gives us $v=\tilde{v}$.

\medskip

2.2. The proof of injectivity of $S_{\eta}$ for $V=\mathbb{R}^{nk}$ is
simpler due to convexity.

\medskip

2.3. From the estimates
$$
|id-DS_{\eta}(v)|_{L_c(\mathbb{R}^{nk},\mathbb{R}^{nk})}=|D_2R(\eta,v)|_{L_c(\mathbb{R}^{nk},\mathbb{R}^{nk})}\le\frac{1}{2}<1\quad\mbox{for}\quad v\in V
$$
we infer that each $DS_{\eta}(v)$, $v\in V$, is an isomorphism. Therefore the Inverse Mapping Theorem applies and yields that $S_{\eta}^{-1}$ is given by continuously differentiable maps on neighbourhoods of the values $y\in S_{\eta}(V)$.

\medskip

3. On assertion (iii). For every $\eta\in W$ and $y=S_{\eta}(v)$ with $v\in V$ we have that $v=S_{\eta}^{-1}(y)$ satisfies
$$
y-(v-R(\eta,v))=0,
$$
or equivalently,
$$
F(\eta,y,v)=0
$$
for the continuously differentiable map $F:(\cup_{\eta\in W}\{\eta\}\times S_{\eta}(V))\times V\to\mathbb{R}^{nk}$ given by
$$
F(\eta,y,v)=y-(v-R(\eta,v)).
$$
Because of the estimates
$$
|D_3F(\eta,y,v)-id|_{L_c(\mathbb{R}^{nk},\mathbb{R}^{nk})}=|D_2R(\eta,v)|_{L_c(\mathbb{R}^{nk},\mathbb{R}^{nk})}\le\frac{1}{2}<1
$$
for $\eta\in W,y\in S_{\eta}(V),$ and $v\in V,$ each map $D_3F(\eta,y,v)$ is an isomorphism. The Implicit Function Theorem applies and yields that the map
$$
\cup_{\eta\in W}\{\eta\}\times S_{\eta}(V)\ni(\eta,y)\mapsto S_{\eta}^{-1}(y)\in\mathbb{R}^{nk}
$$
is locally given by continuously differentiable maps. $\Box$

\medskip

Using Proposition 3.2 (i) and continuity we obtain that the set
$$
{\mathcal O}=\{\chi\in C^1_n:L\chi\in W\,\,\mbox{and}\,\,\widehat{\chi}\in S_{\eta}(V)\,\,\mbox{for}\,\,\eta=L\chi\} 
$$
is open. The map 
$$
B:{\mathcal O}\to C^1_n
$$
given by 
$$
B(\chi)=\chi+g(v)\cdot H(v)\quad\mbox{for}\quad v=S_{\eta}^{-1}(y),\,\,y=\widehat{\chi},\,\,\eta=L\chi
$$ 
is continuously differentiable. Analogously to Eq. (4) we have 
$$
LB(\chi)=L\chi\quad\mbox{for every}\quad\chi\in{\mathcal O}.
$$

\begin{proposition}
$A(U)\subset{\mathcal O}$ and $B(A(\phi))=\phi$ for all $\phi\in U$.
\end{proposition}

{\bf Proof.} Let $\phi\in U$ be given and set $\chi=A(\phi)$. Then $L\chi=L\phi\in W$, and $y=\widehat{\chi}$ satisfies
$y=S_{\eta}(v)$ for $v=\widehat{\phi}$ and $\eta=L\phi$. Hence $\widehat{\chi}\in S_{\eta}(V)$. It follows that $A(\phi)=\chi\in{\mathcal O}$, and we obtain $A(U)\subset{\mathcal O}$. With $\phi,\chi,y,v,\eta$ as before, $v=S_{\eta}^{-1}(y)$, and thereby,
$$
B(A(\phi))=B(\chi)=\chi+g(v)\cdot H(v)=[\phi-g(v)\cdot H(v)]+g(v)\cdot H(v)=\phi.\quad\Box
$$

\begin{proposition}
$B({\mathcal O})\subset U$ and $A(B(\chi))=\chi$ for all $\chi\in{\mathcal O}$.
\end{proposition}

{\bf Proof.} 1. In order to obtain $B({\mathcal O})\subset U$ let $\chi\in{\mathcal O}$ be given and set $\phi=B(\chi)$, $\eta=L\chi$, $y=\widehat{\chi}$. As $L\phi=L\chi\in W$ we may consider $\widehat{\phi}\in\mathbb{R}^{nk}$. In order to show $\phi\in U$ we need to
verify $\widehat{\phi}\in V$. From $\chi\in{\mathcal O}$ we have $y=S_{\eta}(v)$ for some $v\in V$. For every $\iota=(\kappa-1)n+\nu$ with $\kappa\in\{1,\ldots,k\}$ and $\nu\in\{1,\ldots,n\}$,
\begin{eqnarray}
\widehat{\phi}_{\iota} & = & \phi_{\nu}(-d_{\kappa}(L\phi))\nonumber\\
& = & \chi_{\nu}(-d_{\kappa}(L\phi))+g_{\nu}(S_{\eta}^{-1}(y))(H_{\nu}(S_{\eta}^{-1}(y)))(-d_{\kappa}(L\phi))\nonumber\\
& = & \chi_{\nu}(-d_{\kappa}(L\chi))+g_{\nu}(v)(H_{\nu}(v))(-d_{\kappa}(L\chi))\nonumber\\
& = & y_{\iota}+g_{\nu}(v)(H_{\nu}(v))(-d_{\kappa}(\eta))\nonumber\\
& = & y_{\iota}+R_{\iota}(\eta,v)\nonumber\\
& = & S_{\iota}(\eta,v)+R_{\iota}(\eta,v)=v_{\iota},\nonumber
\end{eqnarray}
which yields $\widehat{\phi}=v\in V$.

2. For $\chi,\phi,\eta,y,v$ as in Part 1 of the proof we saw that
$\widehat{\phi}=v$. It follows that
\begin{eqnarray}
A(B(\chi)) & = & A(\phi)=\phi-g(\widehat{\phi})\cdot H(\widehat{\phi})\nonumber\\
& = & [\chi+g(S_{\eta}^{-1}(y))\cdot H(S_{\eta}^{-1}(y))]-g(v)\cdot H(v)\nonumber\\
& = & \chi\quad\mbox{(with}\quad v=S_{\eta}^{-1}(y)).\quad \Box\nonumber
\end{eqnarray}

\medskip

\begin{theorem}
The map $A$ is an almost graph diffeomorphism with respect to
the continuously differentiable submanifold $X_f$, to
the closed subspace $X_0$ of codimension $n$ in $C^1_n$, and to
the open subset $U\supset X_f$ of $C^1_n$, and the solution manifold $X_f$ is an almost graph over $X_0$.	
\end{theorem}

\end{document}